\documentclass[twoside,notitlepage,11pt]{article}

\usepackage{amssymb}
\usepackage[leqno]{amsmath}
\usepackage{amsfonts}
\usepackage{amsopn}
\usepackage{amstext}
\usepackage{amsthm}

\usepackage[backref,colorlinks]{hyperref}

\providecommand{\cal}{\mathcal}
\renewcommand{\Bbb}{\mathbb}

\newenvironment{pf}{\begin{proof}}{\end{proof}}

\newcommand{\Bee}{{\cal{B}}}
\newcommand{\Cee}{{\cal{C}}}

\newcommand{\Ef}{{\cal{F}}}

\newcommand{\Haa}{{\cal{H}}}

\newcommand{\Pee}{{\cal{P}}}

\newcommand{\Nat}{{\Bbb{N}}}
\newcommand{\Qyu}{{\Bbb{Q}}}
\newcommand{\Err}{{\Bbb{R}}}

\newcommand{\lam}{{\lambda}}
\newcommand{\al}{\alpha}

\newcommand{\sig}{\sigma}
\newcommand{\eps}{\varepsilon}
\renewcommand{\phi}{\varphi}
\renewcommand{\rho}{\varrho}

\newcommand{\rest}{\restriction}

\newcommand{\ntr}{n\in\omega}

\newcommand{\loe}{\leqslant}
\newcommand{\goe}{\geqslant}

\newcommand{\subs}{\subseteq}
\newcommand{\sups}{\supseteq}
\newcommand{\nnempty}{\ne\emptyset}

\renewcommand{\iff}{\Longleftrightarrow}

\newcommand{\id}[1]{\operatorname{id}_{#1}}

\newcommand{\dom}{\operatorname{dom}}
\newcommand{\rng}{\operatorname{rng}}

\newcommand{\meet}{\wedge}

\newcommand{\Ht}{\operatorname{ht}}

\newcommand{\concat}{{}^\smallfrown}

\newcommand{\wQyu}{w\Qyu} 
\newcommand{\sigQyu}{\sig\Qyu} 

\newcommand{\cmp}{\circ}
\newcommand{\ciag}[1]{\sett{{#1}_n}{\ntr}}

\newtheorem{tw}{Theorem}[section]
\newtheorem{wn}[tw]{Corollary}
\newtheorem{lm}[tw]{Lemma}
\newtheorem{prop}[tw]{Proposition}

\theoremstyle{definition}

\theoremstyle{remark}

\newcommand{\setof}[2]{\{#1\colon #2\}}

\newcommand{\seq}[1]{\langle #1 \rangle}

\newcommand{\sett}[2]{\{#1\}_{#2}}
\newcommand{\sn}[1]{\{#1\}} 
\newcommand{\dn}[2]{\{#1,#2\}} 
\newcommand{\pair}[2]{{\langle #1, #2 \rangle}} 
\newcommand{\map}[3]{#1\colon #2 \to #3} 
\newcommand{\img}[2]{#1[#2]} 
\newcommand{\inv}[2]{{#1}^{-1}[#2]} 

\newcommand{\power}[1]{\Pee(#1)}
\newcommand{\dpower}[2]{[#1]^{#2}}

\newcommand{\En}{\mathcal N}
\newcommand{\Em}{\mathcal M}

\newcommand{\I}{\ensuremath{\mathcal I}}

\newcommand{\iso}{\cong}

\newcommand{\cont}{\ensuremath{2^{\aleph_0}}}

\newcommand{\Bone}[1]{\ensuremath{\operatorname{B}_1(#1)}}

\newcommand{\suc}[1]{{#1}^+}

\newcommand{\iseg}{\sqsubseteq}

\newcommand{\cee}[1]{{\mathcal C\left(#1\right)}}
\newcommand{\ceezero}[1]{{\mathcal C_0\left(#1\right)}}

\newcommand{\A}[1]{{\alpha}{#1}}
\newcommand{\Af}[2]{{\alpha}_{#2}{#1}}

\title{Finitely fibered Rosenthal compacta and trees}
\author{
{\sc Wies{\l}aw Kubi\'s}
\footnote{Supported in part by the Grant IAA 100 190 901
and by the Institutional Research Plan of the Academy of Sciences of Czech Republic No. AVOZ 101 905 03.}
\\
{\small Institute of Mathematics}\\
{\small Academy of Sciences of the Czech Republic}\\
{\small \v Zitn\'a 25, 115 67 Praha 1}\\
{\small Czech Republic}\\
{\small \it and }\\
{\small Institute of Mathematics}\\
{\small Jan Kochanowski University}\\
{\small Kielce, Poland}
\and
{\sc An\'\i bal Molt\'o}\\
{\small University of Valencia, Spain}\\
{\small Departamento de An\'alisis Matem\'atico}\\
{\small Facultad de Matem\'aticas}\\
{\small Universidad de Valencia}\\
{\small Dr. Moliner 50, 46100 Burjassot (Valencia)}\\
{\small Spain}
}

\begin{document}
\maketitle

\begin{abstract}
We study some topological properties of trees with the interval topology. In particular, we characterize trees which admit a $2$-fibered compactification and we present two examples of trees whose one-point compactifications are Rosenthal compact with certain renorming properties of their spaces of continuous functions.



\end{abstract}

\tableofcontents

\section{Introduction}

A compact space is {\em $n$-fibered} if it has an at most $n$-to-$1$ continuous map onto a metric space. 
Obvious variations of the above definition give the notions of {\em finitely/metrizably fibered} spaces.
A famous open question, attributed to Fremlin (see \cite{Gru}), asks whether it is consistent with the usual axioms of set theory that every perfectly normal compact is $2$-fibered.

Another motivation for studying metrizably fibered compacta is a remarkable result of Todor\-\v ce\-vi\'c~\cite{To2}: every hereditarily separable Rosenthal compact is $2$-fibered.
It seems to be unknown whether every Rosenthal compact is a continuous image of a finitely fibered (or at least metrizably fibered) compact space.

There are natural, more general than finitely (resp. metrizably) fibered, classes of compacta, called in this paper {\em finitely} ({\em metrizably}) {\em determined} spaces.
They include all finitely (resp. metrizably) fibered compacta, while additionally they are stable under continuous images.
These classes were first formally studied independently by Tkachuk~\cite{Tkachuk} and Tka\v cenko~\cite{Tkacenko}.
It has been proved in \cite[Prop. 2.14]{KOSz} that every compact subspace of $\Err^X$ (with the pointwise topology) consisting of functions with countably many points of discontinuity, where $X$ is a fixed Polish space, is metrizably determined.
This seems to confirm the conjecture that every Rosenthal compact is a continuous image of a metrizably fibered compact.

Finitely fibered compacta are in some sense close to metric spaces, so it is natural to ask whether their spaces of continuous functions have good renorming properties (see Question 6.28 in~\cite{MOTV}).
We show that this may not be the case. Namely, there exists a $2$-fibered scattered Rosenthal compact $K$ for which $\cee K$ fails to have a Kadec renorming. We also show that the existence of a Kadec renorming for $\cee K$, with $K$ Rosenthal, does not imply that $K$ is a continuous image of a $2$-fibered compact.
The remaining question is whether spaces of continuous functions over $2$-fibered compacta have rotund renormings. 
One has to mention a recent positive result in this direction \cite{HMO}: $\cee K$ has a locally uniformly rotund renorming whenever $K$ is a separable Rosenthal compact homeomorphic to a set of real functions on a fixed Polish space, where each of the functions has only countably many points of discontinuity.
Since every hereditarily separable Rosenthal compact is $2$-fibered, one may ask whether this renorming result can be proved for all $2$-fibered Rosenthal compacta.
As we have already mentioned, this is not the case.

Our work is also related to the book~\cite{MOTV} and, in particular to Question 6.28 therein asking, roughly speaking, when $\cee K$, with $K$ metrizably determined, has a locally uniformly rotund renorming. Our examples are in the negative direction.

Our inspiration was the work of Todor\v cevi\'c~\cite{T}, where it is proved that the Alexandrov compactification of the tree $\sig\Qyu$ of all bounded well ordered subsets of the rational numbers is Rosenthal compact.
It was known already after Haydon's work \cite{Haydon}, that the Banach space of continuous functions over this particular compact has a rotund renorming, yet it fails to have a Kadec renorming.

We characterize trees whose one point compactification is $2$-determined. Namely, these are precisely $\Err$-embeddable trees of cardinality at most the continuum.
We also give a criterion for being $3$-determined. 
Our examples come from some constructions based on the tree $\sig\Qyu$. We rely on the results of Haydon~\cite{Haydon}, where several renorming properties of spaces of continuous functions over one-point compactifications of trees have been characterized.

\section{Preliminaries}


By a {\em space} we mean a Hausdorff topological space. Let $n>0$ be a natural number. A compact space $K$ is {\em $n$-fibered} if there is a continuous map $\map fKS$ such that $S$ is a second countable space and $|f^{-1}(y)|\loe n$ for every $y\in S$.
We shall say that $K$ is {\em $n$-determined} if there are a second countable space $S$ and an upper semicontinuous map $\map \Phi S{\dpower Kn}$ such that $K=\bigcup_{s\in S}\Phi(s)$. This class of spaces was denoted by $L\Sigma(\loe n)$ in \cite{KOSz}.
It is well known and not hard to prove (see e.g. \cite{KOSz}) that a regular space $K$ is $n$-determined if there are a cover $\Cee\subs\dpower Kn$ and a countable family of closed sets $\En$ which forms a {\em network} for $\Cee$, i.e. for every $C\in\Cee$ and an open set $V\sups C$ there is $N\in\En$ with $C\subs N \subs V$. 
If $K$ is compact, it is enough to require that for every $C\in \Cee$ there is $\En_C\subs \En$ such that $C=\bigcap \En_C$.
Note that every continuous image of a compact $n$-fibered space is $n$-determined. The converse is false, see \cite{KOSz}.
Obvious modification of the definition of ``$n$-fibered/determined" leads to the notion of a {\em metrizably fibered/determined} space.
Every metrizably fibered compact is first countable. The one-point compactification of a discrete space of cardinality $\aleph_1$ is an example of a $2$-determined compact that is not metrizably fibered.

Various classes and topological properties of metrizably determined spaces were studied earlier in~\cite{Tkacenko, Tkachuk} and later in~\cite{Okunev1, Okunev2, Okunev3}.

In this paper we are interested in $n$-fibered and $n$-determined compact spaces, where $n\loe3$.

A {\em Rosenthal compact} is a compact space homeomorphic to a subspace of the space of Baire class one functions $\Bone{P}$ endowed with the pointwise topology, where $P$ is a Polish space. We shall use the fact that the characteristic function $1_F$ of a set $F\subs P$ is of Baire class one iff $F$ is at the same time $F_\sig$ and $G_\delta$.

\subsection{Alexandrov-type compactifications}

Let $X$ be a locally compact space which is not compact. The well known {\em Alexandrov compactification} of $X$ is the space $\A X=X\cup\sn\infty$, where $\infty\notin X$ and a neighborhood of $\infty$ is of the form $\A X\setminus F$, where $F\subs X$ is compact.
This construction can be naturally generalized to obtain more complex compactifications of $X$. 
Namely, fix a continuous map $\map fXK$, where $K$ is a compact space.
Denote by $\tau_X$ and $\tau_K$ the topologies of $X$ and $K$ respectively.
For technical reasons, assume that $X\cap K=\emptyset$. We claim that there exists a unique compact topology $\tau$ on $X\cup K$ which extends the topologies of $X$ and $K$ and for which the map $\map r{X\cup K}K$, defined by conditions
$$r\rest X=f\quad\text{ and }\quad r\rest K=\id K,$$
is continuous.

Let us first see uniqueness.
The continuity of $r$ implies that $U\cup \inv fU\in \tau$ for every $U\in\tau_K$. Further, $\tau_X\subs \tau$, because $K$ is closed in $X\cup K$. Finally, $(X\cup K)\setminus F\in\tau$ whenever $F\subs X$ is compact. Now, using the local compactness of $X$, it is straight to check that the family
$$\Bee = \tau_X \cup \setof{(U\cup \inv fU)\setminus F}{U\in\tau_K,\; F\subs X\text{ compact\/}}$$
induces a Hausdorff topology on $X\cup K$. By compactness, $\Bee$ must be a basis of $\tau$. Finally, it is an easy exercise to show that $\pair{X\cup K}\tau$ is indeed compact.

The space $\pair{X\cup K}\tau$ will be denoted by $\Af Xf$. Note that $\Af Xf=\A X$ when $\map fX{\sn\infty}$ is the constant map. If $X$ is a discrete space and $\map fXK$ is any one-to-one map into a compact space, then $\Af Xf$ is the well known Alexandrov duplicate of $\img fX$ in $K$.

We shall need the following property of $\Af Xf$.

\begin{lm}\label{weorjq}
Let $X$ be a locally compact non-compact space and let $f$ be a continuous map from $X$ into a compact space $K$. If both $\A X$ and $K$ are Rosenthal compact then so is $\Af Xf$.
\end{lm}

\begin{pf}
Let $P_0$, $P_1$ be disjoint Polish spaces such that $\A X\subs\Bone{P_0}$ and $K\subs\Bone {P_1}$. We may assume that $\infty\in \A X$ corresponds to the constant zero function $0_{P_0}$ in $\Bone {P_0}$.
Let $P$ be the disjoint topological sum of $P_0, P_1$. Then $P$ is a Polish space. We identify $x\in \Bone{P_1}$ with $x\cup 0_{P_{0}}\in \Bone P$. 
By this way $K\subs\Bone P$.
Define $\map jX{\Bone P}$ by setting
$$j(x)\rest P_0 = x \quad\text{ and }\quad j(x)\rest P_1 = f(x)\rest P_1.$$ 
Clearly, $j$ is well defined, one-to-one and $\img jX \cap K=\emptyset$, because $x\ne 0_{P_0}$ for $x\in X$. Using the continuity of $f$, we conclude that $j$ is a homeomorphic embedding.
It suffices to check that $L = \img jX\cup K$ is closed in $[-\infty,+\infty]^P$.
Since there is only one compact topology on $L$ extending $\img jX$ and $K$ for which $f\cmp j^{-1}$ is continuous, this will ensure that $L$ is homeomorphic to $\Af Xf$.

Define $r_1(x)=x\cdot 1_{P_1}$. That is, $r_1(x)\rest P_1=x\rest P_1$ and $r_1(x)\rest P_0=0$.
Clearly, $\map{r_1}{[-\infty,+\infty]^P}{[-\infty,+\infty]^P}$ is a continuous map.

Fix $v\in [-\infty,+\infty]^P \setminus L$.
Let $v_0=v\rest P_0$. If $v_0 = 0_{P_0}$ then $r_1(v)=v\notin K$ and, using the continuity of $r_1$, we can easily find a neighborhood of $v$ disjoint from $L$. So assume $v_0\ne0_{P_0}$.
Now, if $v_0\notin X$ then, using the compactness of $X\cup\sn{0_{P_0}} = \A X$, we again find a neighborhood of $v$ disjoint from $L$.
It remains to consider the case that $v_0 = \in X$. Notice that $f(v_0)\ne r_1(v)$, because otherwise $v=j(v_0)\in L$. 
Using the continuity of both $f$ and $r_1$, find a neighborhood $V$ of $v$ in $[-\infty,+\infty]^P$ such that $f(x\rest P_0)\ne r_1(x)$ whenever $x\in V$ and $x\rest P_0\in X$.
We may further assume that $x\rest P_0\ne0$ whenever $x\in V$. Then $V\cap K=\emptyset$ and $V\cap \img jX=\emptyset$.
\end{pf}

\subsection{Trees}

For detailed information concerning theory and applications of trees we refer the readers to the excellent survey \cite{To1}. Here we give necessary definitions and facts only.

A {\em tree} is a partially ordered set $\pair T<$ with a minimal element denoted by $0$, such that for every $t\in T$ the set $[0,t)=\setof{x\in T}{x < t}$ is well ordered and for every $s,t\in T$ there exists the greatest lower bound $s\meet t$.
The last condition ensures that its {\em interval topology} induced by open sets of the form $(s,t]=\setof{x\in T}{s<x\loe t}$ is Hausdorff. It is clear that this topology is locally compact, namely every interval of the form $[s,t]=\setof{x\in T}{s\loe x\loe t}$ is compact. We denote by $\A T$ the Alexandrov one-point compactification of the tree $T$.
That is, $\A T = T\cup \sn\infty$ and a basic neighborhood of $\infty$ is of the form
$$\A T \setminus \bigcup_{i<n}[0,t_i],$$
where $t_0,\dots, t_{n-1}\in T$, $\ntr$.
It will be convenient to extend the partial order of $T$ onto $\A T$ by setting $t<\infty$ for every $t\in T$.
By this way, $\infty = \sup\Cee$ whenever $\Cee$ is an unbounded chain in $T$.

Given a tree $T$, the order type of $[0,t)$ is called the {\em height} of $t$ in $T$ and denoted by $\Ht_T(t)$ or just $\Ht(t)$. The {\em height} of $T$ is the supremum of all numbers $\Ht_T(t)$, where $t\in T$. The set $\setof{t\in T}{\Ht_T(t)=\al}$ is called the {\em $\al$-th level} of $T$. A {\em branch through} $T$ is a maximal linearly ordered subset of $T$.
Given $t\in T$ we denote by $\suc t$ the set of all {\em immediate successors} of $t$ in $T$, that is, $\suc t=\setof{s\in T}{t<s\land [t,s] = \dn ts}$. We say that $T$ is {\em finitely}/{\em countably branching} if $\suc t$ is finite/countable for each $t\in T$.
A set $A \subs T$ is an {\em antichain} if no two elements of $A$ are comparable, i.e. neither $p<q$ nor $q<p$ holds for distinct $p,q\in A$.

Let $\pair T<$ be a tree and let $\pair X<$ be a linearly ordered set (briefly: a {\em line}). We say that $T$ is {\em $X$-embeddable} if there exists a $<$-preserving function from $T$ into $X$. Note that such a function may not be one-to-one.
We shall be particularly interested in $\Err$-embeddable trees, where $\Err$ denotes the real line.
$\Qyu$-embeddable trees are often called {\em special}. A tree is $\Qyu$-embeddable if and only if it is covered by countably many antichains.

A subset $A$ of a partially ordered set $\pair P<$ is an {\em initial segment} if $(\leftarrow,y]\subs A$ whenever $x\in A$, where $(\leftarrow,x]=\setof{y\in P}{y\loe x}$. If $P$ is a tree, then we say that $A$ is an {\em initial subtree}. 

We shall need the following criterion for continuity of maps defined on trees.

\begin{lm}\label{owiejrie}
Let $\pair T<$ be a tree and let $\map f{\A T}X$ be a continuous map into a topological space $X$. Then $f$ is continuous if and only if it satisfies the following two conditions.
\begin{enumerate}
	\item[(1)] $\lim_{\al<\lam}f(t_\al)=f(t)$, whenever $\lam$ is a regular infinite cardinal and $\sett{t_\al}{\al<\lam}$ is a strictly increasing sequence in $T$ with $\sup_{\al<\lam}t_\al=t$, where $t=\infty$ if $\sett{t_\al}{\al<\lam}$ is unbounded. 
	\item[(2)] $\lim_{n\to\infty}f(t_n)=f(\infty)$, whenever $\ciag t$ is an antichain in $T$.
\end{enumerate}
\end{lm}

\begin{pf}
It is obvious that the above conditions are necessary. Assume $\map f{\A T}X$ satisfies (1) and (2).
It is clear that $f$ is continuous at each point $t\in T$.
Indeed, if for some neighborhood $U$ of $f(t)$, $\inv fU$ were not a neighborhood of $t$, then there would exist an increasing sequence $\sett{s_\al}{\al<\lam}\subs [0,t]\setminus \inv fU$  with $t=\sup_{\al<\lam}t_\al$, where $\lam$ is the cofinality of $[0,t)$. This would contradict (1).

Fix a neighborhood $U$ of $f(\infty)$ and let $B = \A T\setminus \inv fU$. 
If $B$ contains a sequence $\sett{t_\al}{\al<\lam}$ which has no upper bound in $T$ then we have that $\infty=\lim_{\al<\lam}t_\al$ and $U$ witnesses that (1) fails. Thus, every chain in $B$ is bounded in $T$. 
Suppose $B$ contains an infinite antichain $\sett{b_n}{\ntr}$. Then $\infty=\lim_{n\to\infty}\sett{b_n}{\ntr}$ while, on the other hand, $U$ witnesses the failure of (2).

It follows that all antichains in $B$ are finite. It is well known and not hard to prove that a tree with this property can have only finitely many branches, say $S_0,\dots, S_{k-1}$. For each $i<k$ let $a_i$ be an upper bound of $S_i$ in $T$. Then $V=\A T\setminus \bigcup_{i<k}[0,a_i]$ is a neighborhood of $\infty$ such that $\img fV\subs U$.
\end{pf}

We shall consider trees with countable branches only. In this case, condition (1) can be replaced by
\begin{enumerate}
	\item[(1')] $\lim_{n\to\infty}f(t_n) = f(t)$ whenever $\ciag t$ is a strictly increasing sequence in $T$ with $t=\sup_{\ntr}t_n$.
\end{enumerate}

\subsection{Trees of sets}

A {\em tree of sets} indexed by a fixed tree $T$ is a family of sets $\setof{A_t}{t\in T}$ satisfying the following conditions:
\begin{enumerate}
	\item[(3)] $A_t\sups A_s$ whenever $t\loe s$ in $T$.
	\item[(4)] $A_t\cap A_s=\emptyset$ whenever $t,s\in T$ are incomparable.
	\item[(5)] $\bigcap_{s\in C}A_s=A_t$ whenever $C\subs T$ is a chain with $t=\sup C$.
	\item[(6)] $\bigcap_{t\in C}A_t=\emptyset$ whenever $C\subs T$ is an unbounded chain.
\end{enumerate}
Assume further that $X$ is a fixed set such that $A_s\subs X$ for every $s\in T$.
Then $\setof{A_t}{t\in T}$ can be regarded as a topological subspace of the Cantor cube $2^X$, identified with the powerset of $X$ (a set corresponds to its characteristic function).

\begin{prop}\label{pwoierpqiwr}
Let $\setof{A_t}{t\in T}$ be a tree of sets. Then the map $\map f{\A T}{\power X}$, defined by $f(t)=A_t$ for $t\in T$ and $f(\infty)=\emptyset$, is continuous.
\end{prop}

\begin{pf}
It is clear that condition (2) of Lemma~\ref{owiejrie} is satisfied, since any sequence of pairwise disjoint sets converges to the empty set.
Condition (1) is also satisfied, because $\bigcap_{\al<\lam}A_{t_\al}=\lim_{\al<\lam}A_{t_\al}$ whenever $\sett{t_\al}{\al<\lam}$ is increasing.
Thus, the continuity of $f$ follows from Lemma~\ref{owiejrie}.
\end{pf}

A tree of sets $\setof{A_t}{t\in T}$ will be called {\em proper} if $A_t\nnempty$ for every $t\in T$ and $A_t\ne A_s$ whenever $t\ne s$.

\begin{wn}\label{ghuehyr}
Let $\sett{A_t}{t\in T}$ be a proper tree of sets. Then $$\setof{A_t}{t\in T} \cup \sn \emptyset \subs \power{A_0}$$
is homeomorphic to $\A T$.
\end{wn}

Note that every tree is isomorphic to a proper tree of sets. Namely, given a tree $T$, the family $\setof{V_t}{t\in T}$, where $V_t=\setof{s\in T}{t\loe s}$, is a tree of pairwise different nonempty subsets of $T$. The above corollary implies that $\A T$ is homeomorphic to $\setof{V_t}{t\in T}\cup\sn\emptyset$, where $\emptyset$ corresponds to $\infty$.

\subsection{Expanding trees}

We describe a well known operation on a tree that replaces each element by a copy of another fixed tree.

Fix a tree $T$ and another, possibly much smaller, tree $S$. For instance, let $S=2^{<2}$ or $S=2^{<\omega}$.
We would like to insert a copy of $S$ at each node of $T$.
For this aim, for each $t\in T$ choose a tree $\sett{D_s(t)}{s\in S}$ of subsets of $\suc t$ such that $D_0(t)=\suc t$.
Now let $T' = T \cup (T\times (S\setminus \sn0))$, where we declare that 
\begin{enumerate}
	\item[(7)] $t<\pair ts$ for every $s\in S\setminus \sn0$,
	\item[(8)] $\pair ts < \pair t{s'}$ whenever $s<s'$ in $S$ and 
	\item[(9)] $\pair ts < r$ whenever $r\in D_s(t)$.
\end{enumerate}
It is clear that this defines a tree order on $T'$, extending the order of $T$. We shall call it an {\em $S$-expansion} of $T$. This construction of course depends on the choice of $D_s(t)$ for $s\in S$, $t\in T$.

In case where the tree $T$ is already a tree of nonempty sets (so, in particular, the order is reversed inclusion) and $|D_s(t)|>1$ for every $\pair ts \in T\times S$, one can represent the $S$-expansion with respect to $D$ as another tree of sets:
$$T'= T \cup \setof{A_s(t)}{s\in S,\; t\in T}, \quad \text{where } A_s(t) = \bigcup D_s(t).$$
Recall that $D_s(t)$ is a family of pairwise disjoint nonempty sets.

The above construction can of course be generalized in such a way that for each $t\in T$ one adds a different tree $S_t$ above $t$. We shall not need this type of expansions here.

Below we present a sample application of expansions of trees.

\begin{prop}\label{oqwjirki}
Let $T$ be an $\Err$-embeddable tree of cardinality $\loe \cont$ and let $\map fT\Err$ be a $<$-preserving function. Then there exist a countably branching tree $T'$ containing $T$ as a subtree and a $<$-preserving function $\map {f'}{T'}\Err$ such that $f'\rest T = f$.
\end{prop}

\begin{pf}
Given $t\in T$ let $W_n(t) = \setof{x\in t+}{f(x)\goe f(t)+1/n}$ and let $P_n(t) = W_n(t)\setminus W_{n-1}(t)$.
Further, using the fact that $|P_n(t)|\loe\cont$, choose a Cantor tree of sets $\sett{P_s(t,n)}{s\in 2^{<\omega}}$ such that $P_\emptyset(t,n)=P_n(t)$
and $|\bigcap_{s\in C}P_s(t,n)|\loe1$ whenever $C\subs 2^{<\omega}$ is an infinite chain.

Now let $S=\sn\emptyset \cup \omega \cup (\omega\times 2^{<\omega})$ regarded as a tree with the order imposed by conditions
\begin{enumerate}
	\item[(${}^*$)] $\emptyset < n < \pair n\emptyset$ for $\ntr$;
	\item[(${}^*_*$)] $\pair ns \loe \pair mr$ iff $n=m$ and $s\subs r$.
\end{enumerate}
That is, $S$ is obtained by ``joining" countably many copies of the Cantor tree $2^{<\omega}$.

Given $t\in T$, define
$D_\emptyset(t) = \suc t$, $D_n(t) = P_n(t)$ and $D_{\pair ns}(t) = P_s(t,n)$. 

Let $T'$ be the $S$-expansion of $T$ with respect to the partitioning $D$. The tree $T'$ is obviously countably branching. It remains to extend $f$ onto $T\times(S\setminus \sn\emptyset)$.

Fix $t\in T$. Define $f'(t)=f(t)$ and $f'(t,n) = f(t)+1/(2n)$. By assumption, $f(x)\goe f(t)+1/n$ whenever $x\in D_n(t)$, therefore so far defined $f'$ preserves the order.
Given $\ntr$, choose a $<$-preserving function $\map\phi {2^{<\omega}}\Err$ so that $f(t)+1/(2n) < \phi(r) < f(t)+1/n$. Finally, given $s=\pair nr\in S$, define $f'(t,s)=\phi(r)$.
It is clear that $f'$ is $<$-preserving.
\end{pf}

\subsection{$\Err$-embeddable trees and Rosenthal compacta}\label{aspjqpw}

The classical tree $\sigQyu$ is defined to be the set of all bounded well ordered subsets of the rationals $\Qyu$, with the ``being initial segment" order. We write ``$s\iseg t$" for ``$s$ is an initial segment of $t$". Actually, the relation $\iseg$ is defined on all subsets of $\Qyu$: $x\iseg y$ iff $x\subs y$ and $\inf(y\setminus x)\goe \sup(x)$, agreeing that $\sup(\emptyset)=-\infty$ and $\inf(\emptyset)=+\infty$.
The tree $\sigQyu$ is $\Err$-embeddable, which is witnessed by the function $\phi(t) = \sum_{q_n\in t}2^{-n}$, where $\sett{q_n}{\ntr}$ is a one-to-one enumeration of $\Qyu$.
Note that the function $t\mapsto \sup(t)$ is not $<$-preserving.

In some situations it is more convenient to consider the tree $\wQyu$ which consists of all (not necessarily bounded) well ordered subsets of the rationals.
The tree $\wQyu$ is complete in the sense that all chains are bounded from above; in particular every branch of $\wQyu$ has the maximal element---an unbounded well ordered subset of $\Qyu$.

A well known fact (see \cite[Thm. 4]{T}) is that the tree $\sigQyu$ is universal for countably branching $\Err$-embeddable trees. For completeness we include a short proof (different from that of \cite[Thm. 4]{T}).

\begin{prop}\label{wieppr}
Let $T$ be a countably branching tree with a $<$-preserving function $\map fT\Err$. Then there exists a tree embedding $\map\psi T{\sigQyu}$ such that $\img \psi T$ is an initial segment of $\sigQyu$ and $f(t)\goe\sup\psi(t)$ for every $t\in T$.
\end{prop}

\begin{pf}
Denote by $\Haa$ the family of all partial tree embeddings $\map hS{\sigQyu}$ satisfying the above condition, i.e. $f(s)\goe\sup h(s)$ for $s\in S$, such that $S$ is an initial segment of $T$ and given $t\in S$, either $\suc t\subs S$ or $t$ is maximal in $S$.
It is rather clear that the union of a chain of maps in $\Haa$ belongs to $\Haa$. Thus, by the Kuratowski-Zorn Lemma, there exists a maximal element $\psi\in\Haa$. We claim that $\dom(\psi)=T$. Suppose this is not the case and fix a minimal $t\in T\setminus S$, where $S=\dom(\psi)$.
Suppose that $t$ is in the closure of $S$. Fix $s_0<s_1<\dots $ in $S$ with $t=\sup_{\ntr}s_n$ and let $x = \bigcup_{\ntr}\psi(s_n)$. Clearly, $x\in\sigQyu$ because it is well ordered and $\sup(x) =\lim_{n\to\infty}\sup\psi(s_n)\loe \sup_{\ntr}f(s_n) \loe f(t)$. 
Thus $\psi\cup\sn{\pair tx}\in\Haa$, a contradiction.
We conclude that $t$ is not in the closure of $S$. Let $s\in S$ be the immediate predecessor of $t$. Note that $s$ must be a maximal element of $S$ and hence $\suc s\cap S=\emptyset$. We shall extend $\psi$ to $S\cup \suc s$. By assumption, $\suc s$ is countable and hence we may find a one-to-one map $\map j{\suc s}\Qyu$ so that $\sup(s)<j(r)\loe f(r)$ for $r\in \suc s$. Finally, define $\map{\psi'}{S\cup \suc s}{\sigQyu}$ by $\psi'\rest S=\psi$ and $\psi'(r) = \psi(s)\cup\sn{j(r)}$ for $r\in \suc s$.
Clearly $\psi'\in\Haa$, a contradiction.
\end{pf}

The above result says in particular that every countably branching $\Err$-embeddable tree is isomorphic to a subtree of $\sigQyu$. In fact, because of Proposition~\ref{oqwjirki}, every $\Err$-embeddable tree of cardinality $\loe\cont$ can be embedded into $\sigQyu$ (see also \cite[Proof of Thm. 4]{T}).

A theorem of Kurepa~\cite{Kurepa} says that the tree $\sigQyu$ is not $\Qyu$-embeddable. 
A short proof can be sketched as follows (see the proof of Theorem 9.8 in \cite{To1} for a more general argument).

Fix a $<$-preserving function $\map \phi{\sig\Qyu}\Qyu$.
We may assume that $\phi$ is bounded from above.
Define inductively a transfinite sequence in $\sig\Qyu$ by setting $t_0=\emptyset$, $t_{\al+1} = t_\al\cup\sn{\phi(t_\al)}$ and $t_\delta = \bigcup_{\xi<\delta}t_\xi$ whenever $\delta$ is a limit ordinal.
It is straight to see, by transfinite induction, that $\sett{t_\al}{\al<\omega_1}$ is strictly increasing in $\sig\Qyu$.
But this is a contradiction, since $\Qyu$ does not contain a copy of $\omega_1$.

To finish this section, we present a short proof of the result of Todor\v cevi\'c \cite{T} saying that the one-point compactification of the tree $\wQyu$ is Rosenthal compact.

Let $P=\power \Qyu$ be endowed with the Cantor set topology. Given $t\in \wQyu$, define $$A_t = \setof{x\in P}{t\iseg x}.$$
It is straight to check that $A_t$ is closed. Note that $A_t\sups A_s$ whenever $t\iseg s$ and $A_t\cap A_s=\emptyset$ whenever $s$ and $t$ are incomparable.
Given an increasing sequence $\ciag t\subs \wQyu$ with $t=\bigcup_{\ntr}t_n$, we have that $A_t=\bigcap_{\ntr}A_{t_n}=\lim_{n\to\infty}A_{t_n}$.
In other words, $\sett{A_t}{t\in \wQyu}$ is a proper tree of nonempty sets.
By Corollary~\ref{ghuehyr}, $\A{(\wQyu)}$ is homeomorphic to $\setof{A_t}{t\in \wQyu}\cup \sn\emptyset$.
Finally, notice that the characteristic function of each $A_t$ is of the first Baire class
\footnote{The characteristic function of a set is of the first Baire class if and only if this set is at the same time $F_\sig$ and $G_\delta$.}.
This shows that $\A{(\wQyu)}$ is Rosenthal compact.

\section{Main results}

In this section we collect our main results concerning trees, finitely determined compacta, Rosenthal compacta and renorming properties of their spaces of continuous functions.

A result of Todor\v cevi\'c \cite{To2} says that a hereditarily separable Rosenthal compact is $2$-fibered.
It is not known whether every Rosenthal compact is a continuous image of a metrizably fibered compact.
A partial positive result is proved in~\cite[Prop. 2.14]{KOSz}, namely a Rosenthal compact is metrizably determined if it can be represented as a set of real functions on a fixed Polish space, each having countably many points of discontinuity.
Note that for each natural number $n$ the space $X_n=(\A{\omega_1})^n$ is a Rosenthal compact that is $(n+1)$-determined and not $n$-determined \cite[Thm. 4.5]{KOSz}.
That $X_n$ is Rosenthal follows from two facts: $\A{\omega_1}$ is easily seen to be Rosenthal and finite products of Rosenthal compacta are Rosenthal.

\subsection{Trees, 2- and 3-determined compacta}

In this subsection we address the question which trees induce $2$-determined or $3$-determined compacta.

\begin{tw}\label{woetjw}
Let $T$ be a tree. Then $\A T$ is $2$-determined if and only if $T$ is $\Err$-embeddable and $|T|\loe\cont$.
\end{tw}

\begin{pf}
Assume $T$ is $\Err$-embeddable and $|T|\loe\cont$.
We are going to describe a countable closed family $\En$ such that for every $t\in T$, $\dn t\infty=\bigcap\En_t$ for some $\En_t\subs \En$. By compactness, this will show that $\A T=T\cup\sn\infty$ is 2-determined.

Fix a $<$-increasing function $\map hT\Err$. We may assume that $h$ is continuous (define $h'(t)=h(t)$ if $t$ is a successor element and $h'(t)=\sup_{s<t} h(s)$ if $t\in T$ is a limit; then $\map{h'}T{\Err}$ is continuous and $<$-preserving). Let $\I$ denote the collection of all closed rational intervals of $\Err$. Given $A\subs T$ denote by $\min A$ the set of all minimal elements of $A$. For each $I\in\I$ fix a countable collection $\Ef_I$ of subsets of $\min \inv hI$ which separates the points of $\min \inv hI$ (here we use the fact that $|T|\loe \cont$). Given $F\in\Ef_I$ define
$$F'=\sn\infty\cup\setof{t\in \inv hI}{(\exists\; s\in F)\; s\loe t}.$$
Observe that $F'$ is closed. Indeed the set $\inv hI\cup\sn\infty$ is closed and if $t\in \inv hI\setminus F'$ then $(0,t]\cap F'=\emptyset$, by the definition of $F'$.

Define $\En=\setof{F'}{F\in\Ef_I,\;I\in\I}$. Then $\En$ is a countable family consisting of closed subsets of $\A T$. Fix $t\in T$ and fix $s\ne t$. We are going to find $N\in\En$ with $t\in N$ and $s\notin N$.

Assume first that $h(s)\ne h(t)$. Find $I\in \I$ such that $h(t)\in I$ and $h(s)\notin I$. Next, find $F\in \Ef_I$ such that the unique element $s\in\min \inv hI\cap (0,t]$ belongs to $F$. Then $t\in F'\in\En$ and $s\notin F'$.
Assume now that $h(s)=h(t)$. Let $r$ be the maximal element below $s,t$. Find $I\in\I$ such that $s,t\in \inv hI$ and $h(r)\notin I$. Let $s_0,t_0\in\min\inv hI$ be such that $s_0\loe s$ and $t_0\loe t$. Find $F\in \Ef_I$ such that $t_0\in F$ and $s_0\notin F$. Then $t\in F'\in \En$ and $s\notin F'$. This completes the proof of sufficiency.

Now suppose that $\A T$ is $2$-determined. Clearly, $|T|\loe\cont$. 
Further, there exist a cover $\Cee$ of $\A T$ consisting of at most $2$-element sets and a countable family $\En$ of closed sets which is a network for $\Cee$.
Fix the families $\Cee$ and $\En$. We may assume that $\En$ is closed under finite intersections.

Define $$T_0=\setof{t\in T}{\dn t\infty \in \Cee}.$$
We claim that $T\setminus T_0$ has a countable height. For suppose otherwise and choose $t_\al\in T\setminus T_0$ so that $\Ht(t_\al)\goe\al$, $\al<\omega_1$. There is $s_\al\in T$ such that $C_\al=\dn{t_\al}{s_\al}\in \Cee$. Now $U_\al=[0,t_\al]\cup [0,s_\al]$ is a neighborhood of $C_\al$. Find $N_\al\in \En$ with $C_\al\subs N_\al\subs U_\al$. Notice that $$\al\loe\Ht(t_\al)\loe\sup_{x\in N_\al}\Ht(x)=\max\dn{\Ht(t_\al)}{\Ht(s_\al)}<\omega_1.$$
On the other hand, since $\En$ is countable, there should exist $N$ so that $N=N_\al$ for uncountably many $\al<\omega_1$. This is a contradiction.

Thus we may assume that $T=T_0$, since $T\setminus T_0$ is clearly $\Qyu$-embeddable.
Write $\En=\sett{N_n}{\ntr}$. Given $t\in T$, define
$$\phi(t)=\setof{\ntr}{[0,t]\cap N_n=\emptyset}.$$
Clearly, $\phi(t)\subs \phi(s)$ whenever $s\loe t$. Assume $s<t$. Since $\dn t\infty\in\Cee$ and $[0,s]$ is closed and disjoint from $\dn t\infty$, by compactness there exists $N_i\in\En$ such that $\dn t\infty\subs N_i$ and $[0,s]\cap N_i=\emptyset$ (recall that $\En$ is closed under finite intersections). This shows that $\phi(t)\ne\phi(s)$ whenever $s<t$. Finally, setting $h(t)=-\sum_{n\in\phi(t)}2^{-n}$, we see that $T$ is $\Err$-embeddable.
\end{pf}

Before stating the next result, we recall the notion of lexicographic product of lines.
Given two lines $X$ and $Y$, their lexicographic product $X\cdot Y$ is the set $X\times Y$ endowed with the order defined by $\pair{x_0}{y_0}<\pair{x_1}{y_1}$ iff either $x_0<x_1$ or else $x_0=x_1$ and $y_0<y_1$.
In other words, $X\cdot Y$ is built from $X$, replacing each element by a copy of $Y$.
Actually, below we are dealing with very particular lexicographic products: $\Err\cdot\Err$ and $\Err\cdot[0,1]$.
Each of them is embeddable into the other.
The latter one is Dedekind complete, that is, every nonempty bounded set has the supremum in $\Err\cdot[0,1]$.

\begin{tw}\label{pweorpoqwer}
Let $T$ be an $(\Err\cdot\Err)$-embeddable tree of cardinality $\loe\cont$. Then $\A T$ is $3$-determined.
\end{tw}

\begin{pf}
First, replace $\Err\cdot\Err$ by a Dedekind complete line $X:=\Err\cdot[0,1]$ and let $\map pX\Err$ be the projection onto the first coordinate. Note that $p$ is continuous with respect to the interval topology on $X$.
Let $\map fTX$ be $<$-preserving. As in the proof of Theorem~\ref{woetjw}, modify $f$ by setting $f'(t)=\sup_{s<t}f(s)$ when $\Ht_T(t)$ is a limit ordinal and $f'(t)=f(t)$ otherwise. Here we have used the fact that $X$ is complete. 
Thus, we may assume that $h$ is continuous.
Let $\map fT\Err$ be the composition of $h$ and the projection $p$.
Then $f$ is continuous, $\loe$-preserving and for each $\lam\in\Err$ the tree
$$T_\lam = \sn{0_T} \cup f^{-1}(\lam)$$
is $\Err$-embeddable (recall our requirement that every tree must have the minimal element --- that is why we have added $0_T$ into $T_\lam$).
Note that $Y_\lam=f^{-1}(\lam)\cup\sn\infty$ is closed in $\A T$. 
By (the proof of) Theorem~\ref{woetjw}, for each $\lam\in\Err$ there is a family $\sett{F_n(\lam)}{\ntr}$ consisting of closed subsets of $Y_\lam$, whose all maximal intersections are of the form $\dn t\infty$, where $t\in f^{-1}(\lam)$.
Let $F_n$ be the closure in $\A T$ of the union $\bigcup_{\lam\in\Err}F_n(\lam)$.
Let $\En$ be a countable family of closed sets defined in the proof of Theorem~\ref{woetjw}.
Finally, let
$$\Em = \En \cup \setof{F_n}{\ntr}.$$
We claim that all maximal nonempty intersections of elements of $\Em$ are of the form $\{s,t,\infty\}$, where $s\loe t$, $f(s)=f(t)$ and $f(s')<f(s)$ whenever $s'<s$.

It is clear that $\sn\infty = \bigcap\Em$.
Fix $t\in T$ and let $C=\bigcap \setof{M\in\Em}{t\in M}$.
Fix $s\in C$, $s\ne t$.
By the proof of Theorem~\ref{woetjw}, we know that $f(s)=f(t)=f(r)$, where $r=s\meet t$ (otherwise we would be able to separate $s$ from $t$ by an element of $\En$).
Suppose $s$ is a minimal element of $f^{-1}(\lam)$.
Then necessarily $s\loe t$, since otherwise $r<s$ and we would have $f(r)<f(s)$. It remains to show that the case $s\notin\min f^{-1}(\lam)$ is impossible.

For suppose $s\notin \min f^{-1}(\lam)$ and find $\ntr$ such that $t\in F_n(\lam)$ and $s\notin F_n(\lam)$.
Find $u<t$ in $f^{-1}(\lam)$ such that $U=(u,s]$ is disjoint from $F_n(\lam)$. Now observe that $U$ is a neighborhood of $s$ not only in $Y_\lam$ but also in $\A T$.
Consequently, $U\cap F_n(\delta)=\emptyset$ for every $\delta\in \Err$.
In particular, $U\cap F_n=\emptyset$ and hence $s\notin F_n$, $t\in F_n$. This shows that $s\notin C$, a contradiction.
\end{pf}

One could go further this direction and try to characterize trees whose one-point compactifications are $k$-determined for $k>2$. 
We stop here, because our examples will give us only $2$- or $3$-determined compacta.

\begin{tw}
Let $T$ be a tree. The following properties are equivalent.
\begin{enumerate}
	\item[(a)] $T$ is $\Err$-embeddable and $|T|\loe\cont$.
	\item[(b)] $T$ admits a one-to-one continuous map onto a separable metric space.
	\item[(c)] $\A T$ is a continuous image of a $2$-fibered compact.
	\item[(d)] $\A T$ is $2$-determined.
\end{enumerate}
\end{tw}

\begin{pf}
(a)$\implies$(b)
In view of Proposition~\ref{wieppr} and the remarks after it, there exists a tree embedding $\map jT{\wQyu}$. That is, $f$ is a one-to-one map satisfying $f(t)<f(t') \iff t<t'$.
We may modify $f$ so that $\img fT$ becomes closed in $\wQyu$. Indeed, define inductively $\map {f'}T{\wQyu}$ by setting $f'(t)=f(t)$ if the $T$-level of $t$ is a successor ordinal and $f'(t)=\sup_{s<t}f'(s)$ if the $T$-level of $t$ is a limit ordinal. 

Thus, we may assume without loss of generality that $T$ is a closed subtree of $\wQyu$. Let $\map hT{\power\Qyu}$ be the inclusion map. Clearly, $h$ is one-to-one and since $T$ is closed in $\wQyu$, $h$ is continuous.

(b)$\implies$(c)
Let $\map hTK$ be as in (b). Enlarging $K$ if necessary, we may assume that $K$ is a compact metric space. Then $\Af Th$ is a $2$-fibered compact that maps onto $\A T$.

(c)$\implies$(d) This is trivial.

(d)$\implies$(a) Clearly, (d) implies that $|T|\loe\cont$. The fact that $T$ is $\Err$-embeddable is included in Theorem~\ref{woetjw}.
\end{pf}

\begin{wn}\label{swefafqwre}
There exists a $2$-fibered Rosenthal compactification $K_0$ of $\wQyu$ whose remainder is homeomorphic to the Cantor set.
\end{wn}

\begin{pf}
By the implication (a)$\implies$(b) of the above theorem, there is a one-to-one continuous map $\map f{\wQyu}K$ such that $K$ is a compact metric space---in fact, after learning its proof we know that $K=\power \Qyu$ and $f$ is the inclusion map.
Let $K_0=\Af {(\wQyu)}f$. By the result of Todor\v cevi\'c~\cite{T}, $\A {(\wQyu)}$ is Rosenthal compact, therefore by Lemma \ref{weorjq}, $K_0$ is Rosenthal compact too.
Finally, $K_0$ is $2$-fibered, because the canonical retraction of $K_0$ onto $\power \Qyu$ is $2$-to-$1$.
\end{pf}

The above arguments together with Proposition~\ref{wieppr} actually show that every $\Err$-embeddable tree of cardinality $\loe\cont$ has a $2$-fibered compactification which is Rosenthal compact.

\subsection{Two examples}

We shall present two trees $T_1$ and $T_2$, whose one-point compactifications $K_1$ and $K_2$ are $3$-determined Rosenthal compacta and their spaces of continuous functions have certain renorming properties. 

Given a tree $T$, denote by $\ceezero T$ the Banach space of all continuous functions $\map fT\Err$ vanishing at infinity,
i.e. such that for every $\eps>0$ there is a compact set $K\subs T$ satisfying $|f(t)|<\eps$ for $t\in T\setminus K$.
The space $\ceezero T$ is naturally identified with the subspace of $\cee {\A T}$ consisting of all $f\in \cee{\A T}$ such that $f(\infty)=0$.

An important work of Haydon~\cite{Haydon} contains several results on renorming properties of spaces of the form $\ceezero T$, where $T$ is a tree.
In particular, it turns out that the existence of a Kadec renorming of $\ceezero T$ (i.e. a renorming such that the weak and the norm topologies coincide on the unit sphere) is equivalent to $\sig$-fragmentability and it is also equivalent to the existence of a $\loe$-preserving function $\map fT\Err$ which has no bad points, where a point $t\in T$ is called {\em good} for $f$ if there are a finite set $F\subs t^+$ and $\eps>0$ such that $f(s) > f(t)+\eps$ for every $s\in \suc t\setminus F$. As one can easily guess, a point is {\em bad} for $f$ if it is not good.

\begin{prop}\label{wefasfhjlhlwq}
Assume $T$ is an $\Err$-embeddable tree and $\ceezero T$ is $\sig$-fragmentable. Then $T$ is $\Qyu$-embeddable.
\end{prop}

\begin{pf} Let $\map fT\Err$ be a $<$-preserving map and let $\map gT\Err$ be a $\loe$-increasing function with no bad points (which exists by Haydon's result \cite[Thm. 6.1]{Haydon}). Then $h=f+g$ is a $<$-preserving map with the property that for every $t\in T$ there is $\eps(t)>0$ such that
$$h(t)+\eps(t)\loe h(s)\qquad\text{ whenever }\quad s>t.$$
Let $T_n$ be the set of all $t\in T$ such that $\eps(t)\goe1/n$. We claim that $\pair{T_n}{<}$ is a tree of height $\loe \omega$. Indeed, if $\setof{t_\al}{\al\loe\omega}$ were strictly increasing in $T_n$ then $h(t_k)\goe h(t_0) + k/n$ for every $k<\omega$ and therefore $h(t_\omega)$ would not be a real number. It follows that each $T_n$ is a countable union of antichains, which shows that $T$ is $\Qyu$-embeddable.
\end{pf}

Let us note that $\ceezero T$ has a Kadec renorming whenever $T$ is a special tree. Indeed, given a $<$-preserving function $\map fT\Qyu$, one can take an order preserving embedding $\map h\Qyu{C}$, where $C\subs\Err$ is a Cantor set, such that $\img h\Qyu$ consists of all points of $C$ which are isolated from the right. Then $hf$ is $<$-increasing with no bad points, therefore by a theorem of Haydon \cite{Haydon}, $\ceezero T$ has a Kadec renorming.

From Corollary~\ref{swefafqwre} we obtain

\begin{wn}\label{wiehruqr}
There exists a $2$-fibered Rosenthal compact whose space of continuous functions is not $\sig$-fragmentable and, in particular, does not have any equivalent Kadec norm.
\end{wn}

We are now going to construct a tree, which results a $3$-determined Rosenthal compact, not a continuous image of any first countable Rosenthal compact.
The tree will be obtained by a certain $S_2$-expansion of $\sigQyu$, where $S_2=\{\emptyset,0,1\}$ is the unique $3$-element binary tree, i.e. $\emptyset<0$, $\emptyset<1$ and $0,1$ are incomparable.
Given $t\in\sigQyu$, choose $r\in\Qyu$ so that $\sup t < r$ and define intervals
\begin{equation}
I_0(t) = [\sup t, r), \quad I_1(t) = [r,+\infty).
\tag{$*$}\label{eqstarghh}
\end{equation}
Additionally, let $I_\emptyset = [\sup t,+\infty)$. Recall that every $u\in \suc t$ is of the form $u = t\cup\sn x$, where $x$ is a rational number from the interval $I_\emptyset(t)$.
Actually, $\sett{I_s(t)}{s\in S_2}$ is a tree of real intervals indexed by $S_2$. Now define
$$D_s(t) = \setof{t\cup\sn x\in \suc t}{x\in I_s(t)}.$$
Finally, let $T_1 = \wQyu \cup T'$, where $T'$ is the $S_2$-expansion of $\sigQyu$ with respect to $D$.

\begin{tw}\label{sefrqkwr}
There exists a scattered Rosenthal compact $K_1$ satisfying the following conditions.
\begin{enumerate}
	\item[(I)] $K_1$ is $3$-determined, not $2$-determined.
	\item[(II)] $\cee{K_1}$ is not $\sig$-fragmentable.
\end{enumerate}
\end{tw}

\begin{pf}
Let $K_1=\A T_1$, where $T_1$ is the above tree. We 
first show that $K_1$ is Rosenthal compact.
Recall that $T_1$ is represented as $\wQyu\cup (\sigQyu\times\dn01)$, where $t<\pair ti<u$ whenever $t\in\sigQyu$, $i\in\dn01$ and $u\in \suc t$.
Define, as in Subsection~\ref{aspjqpw}, $A_t=\setof{x\in\power\Qyu}{t\iseg x}$.
Then $\sett{A_t}{t\in \wQyu}$ is a tree of closed subsets of the Cantor set $\power\Qyu$.

Given $t\in\sigQyu$, $i\in\dn01$, define
$$A_{\pair ti} = \setof{x\in A_t}{\inf(x\setminus t)\in I_i(t)}.$$
Observe that $A_{\pair t1}$ is a closed set and $A_{\pair t0} = A_t \setminus A_{\pair t1}$ is relatively open in $A_t$.
It follows that both sets are at the same time $F_\sig$ and $G_\delta$ in $\power\Qyu$.
In particular, the characteristic functions of the sets $A_t$, $t\in T_1$ are of the first Baire class and $A_t\ne A_r$ if $t\ne r$.
In order to show that $K_1$ is Rosenthal compact, by Proposition~\ref{pwoierpqiwr}, it suffices to check that $\sett{A_t}{t\in T_1}$ is indeed a (proper) tree of sets.
It is clear that $A_t\cap A_r=\emptyset$ whenever $t\ne r$ and $A_t\sups A_r$ whenever $t\loe r$. 
Given a strictly increasing sequence $\ciag u \subs T_1\setminus \wQyu$, one can find $\ciag t \subs \sigQyu$ such that $u_n < t_n < u_{n+1}$ for every $\ntr$. Thus, conditions (5) and (6) in the definition of the tree of sets follow from the fact that $\sett{A_t}{t\in \wQyu}$ is already a tree of sets.

We now show (I).
Let $\map f{\wQyu}\Err$ be strictly $\iseg$-preserving and define $g(t)=\pair{f(t)}0$ for $t\in \wQyu$, $g(t,i) = \pair{f(t)}1$ for $\pair ti\in \sigQyu\times \dn01$. Then $\map g{T_1}{\Err\cdot\dn01}$ is $<$-preserving, showing that $T_1$ is
$(\Err\cdot\dn01)$-embeddable.
By Theorem~\ref{pweorpoqwer}, $K_1$ is $3$-determined. In order to show that $K_1$ is not $2$-determined, by Theorem~\ref{woetjw}, it suffices to show that $T_1$ is not $\Err$-embeddable.
Suppose $\map f{T_1}\Err$ is $<$-preserving and given $t\in \sigQyu$ let $$\eps(t)=\min_{i=0,1}|f(t,i)-f(t)|.$$
Then $\eps(t)>0$, so choose $r(t)\in\Qyu$ satisfying $f(t)<r(t)<f(t)+\eps(t)$.
Then $\map r{\sigQyu}\Qyu$ is $<$-preserving. On the other hand, $\sigQyu$ is not $\Qyu$-embeddable, according to Kurepa's theorem~\cite{Kurepa}. This shows (I).


Finally, suppose that $\cee{K_1}$ is $\sig$-fragmentable. Since $\ceezero{T_1}\iso\cee{K_1}$, by Haydon's result~\cite[Thm. 6.1]{Haydon}, there exists a $\loe$-preserving function $\map {f_1}{T_1}\Err$ with no bad points. Since $S_2$ is a finite tree, the restriction $f\rest \sigQyu$ has no bad points either. But, $\sigQyu$ does not possess such a function, by the proof of Proposition~\ref{wefasfhjlhlwq}.
This contradiction shows (II).
\end{pf}

A modification of the above tree expansion gives another example of a Rosenthal compact, showing that a positive renorming property does not imply $2$-determination.

\begin{tw}\label{thelastone}
There exists a scattered Rosenthal compact $K_2$ satisfying the following conditions.
\begin{enumerate}
	\item[(I)] $K_2$ is $3$-determined, not $2$-determined.
	\item[(II)] $\cee{K_2}$ has an equivalent Kadec norm.
\end{enumerate}
\end{tw}

\begin{pf}
Let $S_\omega=2^{<\omega}$ be the Cantor tree.
We shall construct a tree $T_2 = \wQyu\cup T'$, where $T'$ is an $S_\omega$-expansion of $\sigQyu$ with respect to $D$ defined below.
As before, let $D_\emptyset(t) = \suc t$.
Let $\sett{I_s(t)}{s\in S_\omega}$ be a fixed Cantor tree of subintervals of the interval $[\sup t,+\infty)$, satisfying $I_\emptyset(t) = [\sup t,+\infty)$ and
\begin{equation}
I_{s\concat0}(t) = [a,c), \quad I_{s\concat1}(t) = [c,b),
\tag{$**$}\label{eqstarigru}
\end{equation}
where $a<c<b$ are such that $I_s(t) = [a,b)$.
Here, $\concat$ denotes the usual concatenation of sequences.
Note that $I_s(t)$ is below $I_r(t)$ (i.e. $x<y$ for every $x\in I_s(t)$, $y\in I_r(t)$) whenever $s$ is lexicographically below $r$ in the tree $S_\omega$.
Finally, define
$$D_{s}(t) = \setof{t\cup \sn x \in \suc t}{x\in I_s(t)}.$$
This finishes the definition of the tree $T_2$.

Let $K_2 = \A {T_2}$. Let $\sett{A_t}{t\in \wQyu}$ be as before and define
$$A_{\pair ts} = \setof{x\in A_t}{\inf(x\setminus t) \in I_s(t)}.$$
Notice that $A_{\pair t{s\concat1}}$ is closed in $A_{\pair ts}$ and $A_{\pair t{s\concat0}} = A_{\pair ts} \setminus A_{\pair t{s\concat1}}$.
This implies that all the sets $A_{\pair ts}$ are simultaneously $F_\sig$ and $G_\delta$. 
A similar argument as in the proof of Theorem~\ref{sefrqkwr} shows that $\sett{A_t}{t\in T_2}$ is a proper tree of sets,
which shows that $K_2$ is representable as a space of the first Baire class functions on the Cantor set $\power \Qyu$
(the characteristic functions of $A_t$, $t\in T_2$ plus the constant zero function).

Property (I) is proved like in Theorem~\ref{sefrqkwr}. The only difference is that the tree $T_2$ is $(\Err\cdot \Nat)$-embeddable, not $(\Err\cdot\dn01)$-embeddable.

In order to show (II), notice that the tree $T_2$ is binary, therefore the constant zero function $\map f{T_2}\Err$ has no bad points. By Haydon's theorem~\cite[Thm. 6.1]{Haydon}, we conclude that $\cee{K_2}\iso\ceezero {T_2}$ has a Kadec renorming.
\end{pf}

Let us note in closing that our examples of trees give spaces of continuous functions which do have rotund renormings. This is because of another result of Haydon~\cite[Thm. 5.1]{Haydon}.
Specifically, given a tree $T$, the space $\ceezero T$ has a rotund renorming if and only if there exists an increasing function $\map \rho T\Err$ that is constant on no Cantor subtree of $T$ and such that for every $t\in T$ there is at most one bad point $s\goe t$ with $\rho(s)=\rho(t)$.

Starting with the tree $T=\sig\Qyu$ and making either the $S_2$-expansion or $S_\omega$-expansion, as in the proof of Theorem~\ref{sefrqkwr} and Theorem~\ref{thelastone}, the obvious extension of the function $\rho(t)=\sup t$, $t\in T$, shows that $\cee {K_1}$ and $\cee {K_2}$ have rotund renormings. 
This is because we have chosen particular decompositions of $t^+$, so that from the newly added immediate successors only one of them may become bad
(see conditions (\ref{eqstarghh}) and (\ref{eqstarigru}) respectively).
More precisely, in the case of $T_2$, we need to define
$$\rho(s) = \inf\setof{\rho(t)}{t\in T,\; t\goe s}.$$
Since $T_2$ is binary, there are no bad points and, given $t\in T$, the set $$\setof{s>t}{\rho(s)=\rho(t)} = \{t, \pair t{\seq 0}, \pair t{\seq{00}}, \dots\}$$
is linearly ordered, therefore $\rho$ is not constant on any Cantor subtree of $T_2$.

We do not know whether $K_1$ (or $K_2$) is Rosenthal compact when taking arbitrary decompositions for the tree expansion.

On the other hand, Haydon's work~\cite{Haydon} contains a similar to $T_1$ construction of a tree $\Upsilon$ such that $\ceezero\Upsilon$ fails to have a rotund renorming.
To be more precise, $\Upsilon$ is an $S_2$-expansion of the tree
$$\Gamma = \setof{t\in \omega^{<\omega_1}}{t\text{ is one-to-one and }|\omega\setminus\rng(t)| = \aleph_0}.$$
The order is inclusion or, in other words, extension of functions. The tree $\Gamma$ is easily seen to be $\Err$-embeddable: the function $h(t) = \sum_{n\in\rng(t)}2^{-n}$ is strictly order preserving.
Moreover, it is not hard to prove that $\Gamma$ contains an isomorphic copy of $\sigQyu$. On the other hand, the spaces $\A\Gamma$ and $\A{\sigQyu}$ are not homeomorphic, because the existence of a rotund renorming distinguishes their spaces of continuous functions.

Note that $\A \Gamma$ is a $2$-determined Rosenthal compact.
We do not know whether $\A \Upsilon$ is Rosenthal compact.

\end{document}